\documentclass[a4paper,12pt]{article}

\usepackage{amssymb}

\usepackage{amsmath}
\usepackage{amsfonts}
\usepackage{amscd}
\usepackage[english]{babel}
\title{Dual leftinvariant almost complex structures on the $SU(2)\times SU(2)$.}
\author{Natalia Daurtseva\\ Russia, Kemerovo State University\\
natali0112@ngs.ru}
\date{ }
\begin{document}
\maketitle
\begin{abstract}
There exist non-degenerate 3-form $d\omega_I$, $\omega_I(X,Y)=g(IX,Y)$, for each leftinvariant almost Hermitian structure
$(g,I)$, where $g$ is Killing-Cartan metric on the $M=S^3\times S^3=SU(2)\times SU(2)$. Known \cite{H1},
that arbitrary non-degenerate 3-form on the 6-dimensional manifold, with some additional properties defines the almost complex structure.
Condition for $I$ to define almost complex structure $J_I$ by $d\omega_I$ is obtained.
Properties of $J_I$ are researched.

\end{abstract}

{\bf Introduction.}  Let $M^6$ be a smooth (i.e. $C^{\infty}$), closed, orientable manifold, with differential 3-form $\psi\in\Lambda^3(M)$.
The form $\psi\in\Lambda^3(M)$ defines $\cite{H1}$ endomorphism $K\in\mbox{End}(TM)\otimes\Lambda^6(M)$ by
$$
K(X)=A(i_X\psi\wedge\psi),
$$
where $A:\Lambda^5\longrightarrow TM\otimes\Lambda^6$ is the isomorphism induced by the exterior product
($i_{A(\varphi)}Vol=\varphi$, where $Vol$ is fixed volume form, $A(\varphi)$ is the suitable
vector field here).
Let denote  $\frac16\mbox{tr}K^2$ as $\tau(\psi)$, then $K^2=\mbox{Id}\otimes\tau(\psi)$.

The group $GL(6,\mathbb{R})$ has two open orbits $\mathcal{O}_1$ and $\mathcal{O}_2$ on $\Lambda^3(\mathbb{R}^6)$.
The stabilizer of the forms in the first orbit is $SL(3,\mathbb{C})$. Known \cite{H2}, that $\tau(\psi)<0$
if and only if $\psi\in\mathcal{O}_1$. Thus, if differential 3-form $\psi$ belongs to $\mathcal{O}_1$ at each point,
then it determines the almost complex structure $J$ on $M$:
$$
J=\frac{1}{\kappa}K,
$$
where $\kappa=\sqrt{-\tau(\psi)}$.

Suppose, that now two forms $(\omega,\psi)$, $\omega\in\Lambda^2(M)$ ($\omega\wedge\omega\wedge\omega\neq 0$),
$\psi\in\Lambda^3(M)$ are defined on the $M$. The 3-form $\psi$, with  $\tau(\psi)<0$
provides a reduction to $SL(3,\mathbb{C})$ and 2-form $\omega$ to $Sp(3,\mathbb{R})$.
As the group $SU(3)=Sp(3,\mathbb{R})\bigcap SL(3,\mathbb{C})$,
then two forms $(\omega,\psi)$ with some compatibility conditions provide a reduction to $SU(3)$.

The first of this conditions is:
$$\omega\wedge\psi=0$$
This says that $\omega$ is of type (1,1) with respect to the above almost complex structure $J$.
The second condition is that $\omega(X,JX)$ has to be positive definite form. It gives the Hermitian structure
$(g,\omega,J)$ on $M$, where $g(X,Y)=\omega(X,JY)$.

Further,  if pair $(\omega,\psi)$ gives a reduction to $SU(3)$ and:
$$
\left\{
\begin{array}{l}
\psi=3d\omega;\\
d\phi=-2\mu\omega\wedge\omega,\ \mbox{where}\ i_X\psi=i_{JX}\phi,\ \mu\in\mathbb{R}
\end{array}
\right.
$$
then $(M,g,J)\in\mathcal{NK}$ \cite{Bu}.
This approach is used in  \cite{Bu} to construct invariant nearly K\"{a}hler structure on
$S^3\times S^3$.

Let $M=S^3\times S^3$ now. We'll take interest in leftinvariant structures on $S^3\times S^3=SU(2)\times SU(2)$.
In this case all the calculations are reduced to ones on the Lie algebra
$\mathfrak{su}(2)\times\mathfrak{su}(2)$ of Lie group $SU(2)\times SU(2)$. Denote $(e_1,e_2,e_3,e_4,e_5,e_6)$
the standard frame of $\mathfrak{su}(2)\times\mathfrak{su}(2)=\mathbb{R}^3\times\mathbb{R}^3$
($[e_1,e_2]=e_3$, $[e_1,e_3]=-e_2$, $[e_2,e_3]=e_1$, $[e_4,e_5]=e_6$, $[e_4,e_6]=-e_5$, $[e_5,e_6]=e_4$,
$[e_i,e_j]=0$, for other $i,j$; $e_i\in\mathfrak{su}(2)\times\{0\},\ i=1,2,3$;
$e_i\in\{0\}\times\mathfrak{su}(2),\ i=4,5,6$). Let fix the orientation on $M$, which defined by
the selection of vectors
$(e_1,e_2,e_3,e_4,e_5,e_6)$. Consider the space $\mathcal{A}^+$ of all leftinvariant almost complex structures on $M$
that induce the given orientation.

Fix the Riemannian metric $g$ induced by the Killing-Cartan form on $SU(2)\times SU(2)$.
Take the space $\mathcal{AO}_g^+\subset\mathcal{A}^+$ of all leftinvariant
$g$-orthogonal almost complex structures $J$, induced the given orientation on $SU(2)\times SU(2)$.
Known \cite{D1}, that $\mathcal{A}^+$ is bundle over the $\mathcal{AO}^+_g$,
fiber of the bundle over the $J\in\mathcal{AO}^+_g$ is  the space
$\mathcal{A}^+_{\omega_J}=\{I\in\mathcal{A}^+:\ \omega_J(IX,IY)=\omega_J(X,Y),\ \omega_J(X,IX)>0, \forall X,Y\in\mathfrak{su}(2)\times\mathfrak{su}(2)\}$
of almost complex structures, positive associated with non-degenerate form
$\omega_J(X,Y)=g(JX,Y)$.

The leftinvariant 2-form $\omega_I(X,Y)=g(IX,Y)$ corresponds to each almost complex structure
$I\in\mathcal{AO}^+_g$. Known that closed non-degenerate differential 2-form does not exist on $SU(2)\times SU(2)$. Thus
$d\omega_I\neq 0$.
Therefore the map:
$$
I\in\mathcal{AO}^+_g\longrightarrow d\omega_I\longrightarrow K\in\mbox{End}(\mathfrak{su}(2)\times\mathfrak{su}(2))
$$
is defined.
If $\tau(d\omega_I)<0$, then $I\in\mathcal{AO}^+_g$
defines almost complex structure $J_I\in\mathcal{A}^+$.

The almost complex structures $I\in\mathcal{AO}^+_g$, for which $\tau(d\omega_I)<0$ are investigated in the paper.
The properties of corresponding almost complex structures $J_I\in\mathcal{A}^+$ are studied.

{\bf Calculation.} Let $I\in\mathcal{AO}^+_g$. It is defined by skew-symmetric matrix
$$
I=\left(
\begin{array}{cc}
A & B\\
-B^T & C
\end{array}
\right),
$$
where $A=\left(
\begin{array}{ccc}
0 & a_1 & a_2\\
-a_1 & 0 & a_3\\
-a_2 & -a_3 & 0
\end{array}
\right)$, $B=
\left(
\begin{array}{ccc}
b_1 &  b_2 &  b_3\\
b_4 &  b_5 &  b_6\\
b_7 &  b_8 &  b_9
\end{array}
\right)$,
$C=\left(
\begin{array}{ccc}
0 &  c_1 &  c_2\\
-c_1 &  0 &  c_3\\
-c_2 &  -c_3 &  0
\end{array}
\right).
$ in the standard frame.
Parameters $a_i,b_j,c_k$, $i,k=1,2,3$, $j=1,..,6$ are in relations given by condition $I^2=-1$, in particular we have:
$$
\begin{array}{l}
1)\ b_1^2+b_2^2+b_3^2+a_1^2+a_2^2=1.\\
2)\ b_4^2+b_5^2+b_6^2+a_1^2+a_3^2=1.\\
3)\ b_7^2+b_8^2+b_9^2+a_2^2+a_3^2=1.\\
4)\ a_2a_3+b_1b_4+b_2b_5+b_3b_6=0.\\
5)\ -a_1a_3+b_1b_7+b_2b_8+b_3b_9=0.\\
6)\ a_1a_2+b_4b_7+b_5b_8+b_6b_9=0.\\
7)\ b_1^2+b_4^2+b_7^2+c_1^2+c_2^2=1.\\
8)\ b_2^2+b_5^2+b_8^2+c_1^2+c_3^2=1.\\
9)\ b_3^2+b_6^2+b_9^2+c_2^2+c_3^2=1.\\
\end{array}
$$

{\bf Theorem 1.} {\it $\tau(d\omega_I)<0$ if and only if $a_1^2+a_2^2+a_3^2<\frac34$.}

{\bf Proof.} Let use the Maurer-Cartan formulas to calculate $d\omega_I$:
$d\theta^k=-\sum_{i<j}C^k_{ij}\theta^i\wedge\theta^j$, where $C^k_{ij},\ k,i,j=1\dots6$ are the structural constants,
$\{\theta^i, i=1\dots6\}$ the frame of the leftinvariant 1-forms space on $M=SU(2)\times SU(2)$. Let $(e^1,e^2,e^3,e^4,e^5,e^6)$ is
co-frame to $(e_1,e_2,e_3,e_4,e_5,e_6)$. Then $de^1=-e^2\wedge e^3$, $de^2=-e^3\wedge e^1$, $de^3=-e^1\wedge e^2$,
$de^4=-e^5\wedge e^6$, $de^5=-e^6\wedge e^4$, $de^6=-e^4\wedge e^5$. We have
$$
\psi=d\omega_I=b_1(e^{234}-e^{156})+b_2(e^{235}-e^{164})+b_3(e^{236}-e^{145})+b_4(e^{314}-e^{256})+
$$
$$
+b_5(e^{315}-e^{264})+b_6(e^{316}-e^{245})+b_7(e^{124}-e^{356})+b_8(e^{125}-e^{364})+b_9(e^{126}-e^{345})
$$
where $e^{ijk}$ is the 3-form $e^i\wedge e^j\wedge e^k$.
$$
i_{e_1}\psi\wedge\psi=(-b_1^2-b_2^2-b_3^2+b_4^2+b_5^2+b_6^2+b_7^2+b_8^2+b_9^2)e^{23456}+
$$
$$
+2(b_1b_4+b_2b_5+b_3b_6)e^{13456}-2(b_1b_7+b_2b_8+b_3b_9)e^{12456}+
$$
$$
+2(b_6b_8-b_5b_9)e^{12356}+2(b_6b_7-b_4b_9)e^{12346}+2(b_5b_7-b_4b_8)e^{12345}
$$
Then
$$
i_{e_1}\psi\wedge\psi=(1-2a_3^2)e^{23456}-2a_2a_3e^{13456}-2a_1a_3e^{12456}+2(b_6b_8-b_5b_9)e^{12356}+
$$
$$
+2(b_6b_7-b_4b_9)e^{12346}+2(b_5b_7-b_4b_8)e^{12345}
$$

As $i_{K(X)}Vol=i_X\psi\wedge\psi$, and $i_Y Vol=Y^1e^{23456}-Y^2e^{13456}+Y^3e^{12456}-Y^4e^{12356}+
Y^5e^{12346}-Y^6e^{12345}$, where $Vol=e^1\wedge e^2\wedge e^3\wedge e^4\wedge e^5\wedge e^6$, then matrix of
$\frac12 K$ in the above basis is
$$
\left(
\begin{array}{cccccc}
\frac12-a_3^2 & a_2a_3 & -a_1a_3 & b_8b_6-b_9b_5 & b_9b_4-b_7b_6 & b_7b_5-b_8b_4\\
a_2a_3 & \frac12-a_2^2 & a_1a_2 & b_9b_2-b_8b_3 & b_7b_3-b_1b_9 & b_8b_1-b_7b_2\\
-a_1a_3 & a_1a_2 & \frac12-a_1^2 & b_5b_3-b_6b_2 & b_6b_1-b_4b_3 & b_4b_2-b_5b_1\\
b_9b_5-b_8b_6 & b_8b_3-b_9b_2 & b_2b_6-b_3b_5 & -\frac12+c_3^2 & -c_2c_3 & c_1c_3\\
b_7b_6-b_9b_4 & b_1b_9-b_3b_7 & b_4b_3-b_6b_1 & -c_2c_3 & -\frac12+c_2^2 & -c_1c_2\\
b_8b_4-b_7b_5 & b_2b_7-b_8b_1 & b_5b_1-b_4b_2 & c_1c_3 & -c_1c_2 & -\frac12+c_1^2
\end{array}
\right)
$$
We can calculate $(K^2)_{11}$:
$$
(K^2)_{11}=(1-2a_3^2)^2+4a_2^2a_3^2+4a_1^2a_3^2-4(b_8b_6-b_9b_5)^2-4(b_9b_4-b_7b_6)^2-4(b_7b_5-b_8b_4)^2
$$
Denote
$$
x=a_1^2+a_2^2+a_3^2
$$
Then
$$
(K^2)_{11}=1-4a_3^2+4a_3^2x-4b_8^2b_6^2-4b_9^2b_5^2+8b_8b_6b_9b_5-4b_9^2b_4^2-4b_7^2b_6^2+8b_9b_4b_7b_6-
$$
$$
-4b_7^2b_5^2-4b_8^2b_4^2+8b_7b_5b_8b_4=1-4a_3^2(1-x)-4b_4^2(b_7^2+b_8^2+b_9^2)-4b_5^2(b_7^2+b_8^2+b_9^2)-
$$
$$
-4b_6^2(b_7^2+b_8^2+b_9^2)+4b_4^2b_7^2+4b_5^2b_8^2+4b_6^2b_9^2+8b_8b_6b_9b_5+8b_9b_4b_7b_6+8b_7b_5b_8b_4=
$$
$$
=1-4a_3^2(1-x)-4(b_4^2+b_5^2+b_6^2)(b_7^2+b_8^2+b_9^2)+4(b_4b_7+b_5b_8+b_6b_9)^2
$$
Use the condition $I^2=-1$, obtain:
$$
(K^2)_{11}=1-4a_3^2(1-x)-4(1-a_1^2-a_3^2)(1-a_2^2-a_3^2)+4(-a_1a_2)^2=4x-3
$$
As $(K^2)_{11}=(K^2)_{ii},\ \forall i=1..6$, we have
$$
\tau(d\omega_I)=\frac16\mbox{tr}K^2=-3+4(a_1^2+a_2^2+a_3^2)
$$
\begin{flushright}$\Box$\end{flushright}

{\bf Remark 1.} One can write the matrix of $K$ in the frame $(e)$:
$$
K=\left(
\begin{array}{cc}
1-2A^* & -2B^*\\
2B^{*T} & -1+2C^*
\end{array}
\right),
$$
where  $A^*$, $B^*$, $C^*$ - are the algebraic supplement matrices to  $A$, $B$, $C$ respectively.

{\bf Remark 2.} As $a_1^2+a_2^2+a_3^2=c_1^2+c_2^2+c_3^2$ and $\sum_{i=1}^9b_i^2=3-2\sum_{i=1}^3a_i^2$
by $I^2=-1$, then the conditions:

1. $\tau(d\omega_I)<0$,

2. $c_1^2+c_2^2+c_3^2<\frac34$,

3. $a_1^2+a_2^2+a_3^2<\frac34$,

4. $\sum_{i=1}^9b_i^2>\frac32$\\
are equivalent.

Denote the set of all $I\in\mathcal{AO}^+_g$, with
$\tau(d\omega_I)<0$ as $\mathcal{AO}_-$. Now we can define map $\alpha:\mathcal{AO}_-\longrightarrow\mathcal{A}^+$,
$\alpha(I)=J_I=\frac{1}{\sqrt{-\tau(d\omega_I)}}K$ for each $I\in\mathcal{AO}_-$.

{\bf Remark 3.} $\tau(d\omega_J)=1>0$ for integrable almost complex structure
$J=\left(
\begin{array}{cc}
A & B\\
-B & A
\end{array}
\right),$
$A=\left(
\begin{array}{ccc}
0 & 0 & 0\\
0 & 0 & -1\\
0 & 1 & 0
\end{array}
\right),$
$B=\left(
\begin{array}{ccc}
-1 & 0 & 0\\
0 & 0 & 0\\
0 & 0 & 0
\end{array}
\right)$. For
$I_0=\left(
\begin{array}{cc}
0 & -E\\
E & 0
\end{array}
\right)$
we have $\tau(d\omega_I)=-3$, then $J_{I_0}=\frac{1}{\sqrt{3}}\left(
\begin{array}{cc}
1 & -2\\
2 & -1
\end{array}\right)$ with metric $g_{\mathcal{NK}}=\frac{1}{\sqrt{3}}\left(
\begin{array}{cc}
2 & -1\\
-1 & 2
\end{array}\right)$ gives nearly K\"{a}hler structure.

{\bf Lemma 1. }{\it $\det B<0$ for all almost complex structures $I\in\mathcal{AO}_-$.}

{\bf Proof.} Assume that the almost complex structure $I\in\mathcal{AO}_-$,
such that $\det B=0$ exists. Then we can find non-zero vector $v\in\mathbb{R}^3$, for which $B^Tv=0$.
Let $V\in\mathfrak{su}(2)\times\mathfrak{su}(2)$, such that its projection on the first factor is $v$ ($\pi_1(V)=v$) and
its projection on the second factor is zero ($\pi_2(V)=0$), i.e. $V^T=(v^T,0)$. Then
$$
(I^2V)^T=(v^TA^2,v^TAB)=(-v^T,0)\ \Leftrightarrow\
\left\{
\begin{array}{l}
A^2v=-v;\\
B^TAv=0.
\end{array}
\right.
$$
Condition $A^2v=-v$ is equivalent to $\det (A^2+1)=0$, i.e. $x=1$, this one
contradicts to $x<3/4$.

As $\det B\neq 0$, then vectors $(e_1,e_2,e_3)$ are $I$-linearly non-dependent. Then
$$
(e_1,e_2,e_3,Ie_1,Ie_2,Ie_3)=(e_1,e_2,e_3,e_4,e_5,e_6)\left(
\begin{array}{cc}
E & A\\
0 & -B^T
\end{array}
\right)
$$
I.e. $I$ keeps the orientation \cite{K} if $\det B<0$.
\begin{flushright}$\Box$\end{flushright}

{\bf Lemma 2.} $J_I\in\mathcal{A}^+$.

{\bf Proof.} As $\det B\neq 0$, then $(e_1,e_2,e_3)$ are $J_I$-linearly non-dependent.
The determinant of amplication matrix from the frame $(e_1,e_2,e_3,J_Ie_1,J_Ie_2,J_Ie_3)$ to the standard one is
$$
\det\left(
\begin{array}{cc}
E & 1-2A^*\\
0 & 2B^{*T}
\end{array}
\right)=8(\det B)^2>0.
$$
\begin{flushright}$\Box$\end{flushright}

{\bf Theorem 2. }{\it $J_I\in\mathcal{A}_{\omega_J}^+$, where\\
$$\omega_J=
\frac{2}{\sqrt(1-\tau)}\left(
\begin{array}{cc}
0 & (1+yA^*)B^*\\
-(1+yC^*)B^{*T} & 0
\end{array}\right),\
y=\frac{1-\sqrt{1-x}}{x\sqrt{1-x}}$$}

{\bf Proof.} Find the projection $J_I\in\mathcal{A}^+$ on the base $\mathcal{AO}_g^+$ of bundle $\mathcal{A}^+$,
using results of \cite{D1}.
The skew-symmetric form $\omega(X,Y)=\frac12(g(J_IX,Y)-g(X,J_IY))$ corresponds to structure $J_I$. The matrix of form $\omega$
in the standard frame:
$$
\omega=\frac{2}{\sqrt{-\tau}}\left(
\begin{array}{cc}
0 & B^*\\
-B^{*T} & 0
\end{array}
\right)
$$
Corresponding skew-symmetric operator $D$, such that $\omega(X,Y)=g(DX,Y)$ is:
$$
D=\frac{2}{\sqrt{-\tau}}\left(
\begin{array}{cc}
0 & -B^*\\
B^{*T} & 0
\end{array}
\right)
$$
Then
$$
(-D^2)^{-\frac12}=\left(
\begin{array}{cc}
(-\frac{4B^*B^{*T}}{\tau})^{-1/2} & 0\\
0 & (-\frac{4B^{*T}B^*}{\tau})^{-1/2}
\end{array}
\right)
$$
Find $(B^*B^{*T})^{-1/2}$ and $(B^{*T}B^*)^{-1/2}$. Use condition
$J_I^2=-1$, it is equivalent to system of equations:
$$
\left\{
\begin{array}{l}
(1-2A^*)^2-4B^*B^{*T}=\tau E,\\
B^*C^*=A^*B^*,\\
(1-2C^*)^2-4B^{*T}B^*=\tau E.
\end{array}
\right.
$$
Then $\left(-\frac{4B^*B^{*T}}{\tau}\right)^{-\frac12}=\left(1+\frac{(1-2A^*)^2}{-\tau}\right)^{-\frac12}$ and
$\left(-\frac{4B^{*T}B^*}{\tau}\right)^{-\frac12}=\left(1+\frac{(1-2C^*)^2}{-\tau}\right)^{-\frac12}$.

By direct calculations one can show that $A^{*2}=xA^*$. Then
$$
1+\frac{(1-2A^*)^2}{-\tau}=1+\frac{(1+4(x-1)A^*)}{-\tau}=(1+\frac{1}{-\tau})(1-A^*).
$$
$$
\left(1+\frac{(1-2A^*)^2}{-\tau}\right)^{-\frac12}=\frac{1}{\sqrt{1-\frac{1}{\tau}}}(1-A^*)^{-\frac12}=
$$
$$
=\frac{1}{\sqrt{1-\frac{1}{\tau}}}(1+\frac12 A^*+\frac{3}{2^2}A^{*2}+\frac{15}{2^3}A^{*3}+\dots)=
$$
$$
=\frac{1}{\sqrt{1-\frac{1}{\tau}}}(1+\frac12 A^*+\frac{3}{2^2}xA^*+\frac{15}{2^3}x^2A^*+\dots)=
$$
$$
=\frac{1}{x\sqrt{1-\frac{1}{\tau}}}(x-A^*+A^*(1+\frac12x+\frac{3}{2^2}x^2+\frac{15}{2^3}x^3+\dots))=
$$
$$
=\frac{1}{x\sqrt{1-\frac{1}{\tau}}}(x+A^*(1-x)^{-\frac12}-A^*)=\frac{1}{\sqrt{1-\frac{1}{\tau}}}\left(1+\frac{1-\sqrt{1-x}}{x\sqrt{1-x}}A^*\right).
$$
Projection of almost complex structure $J_I$ on $\mathcal{AO}_g^+$ is:
$$\pi(J_I)=(-D^2)^{-\frac12}D=\frac{2}{\sqrt{1-\tau}}\left(
\begin{array}{cc}
0 & -(1+yA^*)B^*\\
(1+yC^*)B^{*T} & 0
\end{array}
\right),$$
therefore
$\omega_J=\frac{2}{\sqrt{1-\tau}}
\left(
\begin{array}{cc}
0 & (1+yA^*)B^*\\
-(1+yC^*)B^{*T} & 0
\end{array}\right).
$
\begin{flushright}$\Box$\end{flushright}

{\bf Remark 4.} As $\pi(J_I)\in\mathcal{AO}_g^+$, then matrix $X=\frac{2}{\sqrt{1-\tau}}(1+yA^*)B^*\in SO(3)$.
So the matrix of form $\omega$ is canonical $\left(\begin{array}{cc}
0 & E\\
-E & 0
\end{array}\right)$ in the frame $(u_1,u_2,u_3,u_4,u_5,u_6)=(e_1,e_2,e_3,Xe_4,Xe_5,Xe_6)$. The group
$SO(3)\times SO(3)$ acts on $M=S^3\times S^3=SU(2)\times SU(2)$ with canonical metric $g$ by an isometry \cite{Bu}.
We have $[u_4,u_5]=u_6$, $[u_4,u_6]=-u_5$, $[u_5,u_6]=u_4$,
$u_4,u_5,u_6\in\{ 0\}\times\mathfrak{su}(2)$. Therefore frames $(e)$ and $(u)$ are equivalent \cite{Bu}.
One can repeat all the above calculations for the almost complex structure $I\in\mathcal{AO}^+_g$
in the new frame $(u)$, and get $\pi(J_I)=I_0$, where $I_0\in\mathcal{AO}_g^+$: $I_0u_4=u_1$, $I_0u_5=u_2$,
$I_0u_6=u_3$.

{\bf Lemma 3.}{\it The form $\omega_J$ and metric $g_I$
$$
\omega_J=\left(
\begin{array}{cc}
0 & E\\
-E & 0
\end{array}
\right)\qquad
g_I=\left(
\begin{array}{cc}
2t\left(1-\frac{A^*}{1+t}\right) & -1+2A^*\\
-1+2A^* & 2t\left(1-\frac{A^*}{1+t}\right)
\end{array}
\right),
$$
where $t=\sqrt{1-x}$ in the frame $(u_1,u_2,u_3,u_4,u_5,u_6)$.}

{\bf Corollary.} {\it Map
$$
I\in\mathcal{AO}_-\longrightarrow^{\alpha} J_I\longrightarrow^{\pi}\pi(J_I)\longrightarrow^{\alpha} J_{\pi(J_I)}
$$
gives the almost complex structure, for which $(g_{\mathcal{NK}},J_{\pi(J_I)})\in\mathcal{NK}$.}

{\bf Remark 5.} As each almost complex structure $J_I$ for $I\in\mathcal{AO}_-$ is in $\mathcal{A}_{\omega_J}^+$, then
we have almost Hermitian structure $(g_I,J_I)$, where $g_I(X,Y)=\omega_J(X,J_IY)$ for all $I\in\mathcal{AO}_-$.

{\bf Theorem 3.}{\it
Ricci curvature of metric $g_I$ is
$$Ric_{11}=Ric_{22}=-\frac{8t^4-16t^3-10t^2+10t+3}{2(4t^2-1)^2};
Ric_{33}=\frac{4t^2+1}{2(4t^2-1)^2};
$$
$$
Ric_{44}=Ric_{55}=-\frac{8t^4-16t^3+6t^2-2t-1}{2(4t^2-1)^2};
Ric_{66}=\frac{-3+16t^4-8t^2}{2(4t^2-1)^2};
$$
$$
Ric_{ij}=0,\mbox{ if }i\neq j
$$
in the orthonormal frame for $g_I$. The scalar curvature $s=-\frac{8t^4-32t^3-2t^2+8t+3}{(4t^2-1)^2}$.
}

{\bf Proof.}
In case $a_1^2+a_3^2\neq 0$ proper values and proper vectors of matrix $g_I$ in the frame $(u)$ are:
$$
\begin{array}{ll}
\lambda_1=2t-1: & v_1=\frac{1}{\sqrt{2(a_1^2+a_3^2)(2t-1)}}(-a_1,0,a_3,-a_1,0,a_3);\\
 &  v_2=\frac{1}{\sqrt{2(a_1^2+a_3^2)(1-t^2)(2t-1)}}(a_2a_3,a_1^2+a_3^2,a_1a_2,a_2a_3,a_1^2+a_3^2,a_1a_2);\\
\lambda_2=1: & v_3=\frac{1}{\sqrt{2(1-t^2)}}(-a_3,a_2,-a_1,-a_3,a_2,-a_1);\\
\lambda_3=2t+1: & v_4=\frac{1}{\sqrt{2(a_1^2+a_3^2)(2t+1)}}(-a_1,0,a_3,a_1,0,-a_3);\\
 & v_5=\frac{1}{\sqrt{2(a_1^2+a_3^2)(2t+1)}}(a_2a_3,a_1^2+a_3^2,a_1a_2,-a_2a_3,-a_1^2-a_3^2,-a_1a_2);\\
 \lambda_4=4t^2-1: & v_6=\frac{1}{\sqrt{2(1-t^2)(4t^2-1)}}(-a_3,a_2,-a_1,a_3,-a_2,a_1)
\end{array}
$$
In case $a_1^2+a_3^2=0,\ a_2\neq 0$ proper values and proper vectors of matrix $g_I$ in the frame $(u)$ are:
$$
\begin{array}{ll}
\lambda_1=2t-1: & v_1=\frac{1}{\sqrt{2(2t-1)}}(1,0,0,1,0,0);\\
 &  v_2=\frac{1}{\sqrt{2(2t-1)}}(0,0,1,0,0,1);\\
\lambda_2=1: & v_3=\frac{1}{\sqrt{2(1-t^2)}}(0,a_2,0,0,a_2,0);\\
\lambda_3=2t+1: & v_4=\frac{1}{\sqrt{2(2t+1)}}(1,0,0,-1,0,0);\\
 & v_5=\frac{1}{\sqrt{2(2t+1)}}(0,0,1,0,0,-1);\\
 \lambda_4=4t^2-1: & v_6=\frac{1}{\sqrt{2(1-t^2)(4t^2-1)}}(0,a_2,0,0,-a_2,0)
\end{array}
$$
In case $a_1^2+a_2^2+a_3^2=0$ we have $t=1$ and
$$
\begin{array}{ll}
\lambda_1=1: & v_1=\frac{1}{\sqrt{2}}(1,0,0,1,0,0);\\
 &  v_2=\frac{1}{\sqrt{2}}(0,1,0,0,1,0);\\
 & v_3=\frac{1}{\sqrt{2}}(0,0,1,0,0,1);\\
\lambda_2=3: & v_4=\frac{1}{\sqrt{6}}(1,0,0,-1,0,0);\\
 & v_5=\frac{1}{\sqrt{6}}(0,1,0,0,-1,0);\\
 & v_6=\frac{1}{\sqrt{6}}(0,0,1,0,0,-1)
\end{array}
$$
The direct calculations by formulas from \cite{Be} gives tensor Ricci in the frame $(v)=(v_1,v_2,v_3,v_4,v_5,v_6)$.

\end{document}